\nonstopmode
\documentclass[]{amsart}
\usepackage{amssymb}
\usepackage{graphicx}
\usepackage{epsfig}
\usepackage[all]{xy}
\usepackage[dvipsnames]{xcolor}
\usepackage[hyperindex=true,bookmarks=true,bookmarksnumbered=true,debug=true]{hyperref}

\allowdisplaybreaks

\def\undersetbrace#1\to#2{\underbrace{#2}_{#1}}
\def\oversetbrace#1\to#2{\overbrace{#2}^{#1}}
\def\AMSunderset#1\to#2{\underset{#1}{#2}}
\def\AMSoverset#1\to#2{\overset{#1}{#2}}

\newcommand{\nmb}[2]{\ifx!#1{\ref{nmb:#2}}%
\else\if.#1{\label{nmb:#2}}%
\else\if0#1{\label{nmb:#2}}%
\else{{#2}}%
\fi\fi\fi}
\swapnumbers
\newtheorem{proposition}[subsection]{Proposition}
\newtheorem*{proposition*}{Proposition}
\newtheorem{theorem}[subsection]{Theorem}
\newtheorem*{theorem*}{Theorem}
\newtheorem{lemma}[subsection]{Lemma}
\newtheorem*{lemma*}{Lemma}
\newtheorem{corollary}[subsection]{Corollary}
\newtheorem*{corollary*}{Corollary}
\newtheorem{remark}[subsection]{Remark}
\newtheorem*{remark*}{Remark}

\def\ign#1{}             
\def\o{\circ}
\def\X{\mathfrak X}

\def\be{\beta}

\def\Ga{\Gamma}

\def\x{\times}
\def\p{\partial}
\let\on=\operatorname

\def\AMSonly#1{}

\newcommand{\todopm}[1]
  {\vspace{5 mm}\par \noindent \marginpar{\LARGE PM} \framebox{\begin
  {minipage}[c]{0.95 \textwidth} \tt\color{WildStrawberry} #1
\end{minipage}}\vspace{5 mm}\par}
\newcommand{\todopr}[1]
  {\vspace{5 mm}\par \noindent \marginpar{\LARGE PR} \framebox{\begin
  {minipage}[c]{0.95 \textwidth} \tt\color{blue} #1
\end{minipage}}\vspace{5 mm}\par}

\begin{document}
\title[]
{Poisson bivectors on infinite dimensional manifolds
}
\author{Peter W. Michor}
\address{
Peter W. Michor:
Fakult\"at f\"ur Mathematik, Universit\"at Wien,
Oskar-Morgenstern-Platz 1, A-1090 Wien, Austria.
}
\email{Peter.Michor@univie.ac.at}
\author{Praful Rahangdale}
\address{
Praful Rahangdale:
 Institut f\"ur Mathematik, Universit\"at Paderborn,
Warburger Str. 100, 33098 Paderborn,Germany}
\email{praful@math.uni-paderborn.de}
\date{\today}

\begin{abstract} 
We show that, on a smoothly paracompact convenient manifold $M$ modeled on a convenient space with the bornological approximation property,  the dual map of a Poisson bracket factors as a smooth section of the vector bundle $L_{\text{skew}}^2(T^*M,\mathbb R)$.
\end{abstract}
\def\LaTeXonly{}

\maketitle

\section{Introduction} \nmb0{1}
For a finite-dimensional smooth manifold $M$, it is well known that a Poisson bracket on $C^\infty(M)$ induces a Poisson tensor $\pi\in \Gamma(\wedge^2TM)$ (see \cite{drinfeld83}, \cite{weinstein1983local}). But for a manifold modeled on a locally convex space or on a convenient space it is not automatic (see \cite{tumpach2020banach}, \cite{golinski2024poisson}). In \cite{beltictua2018queer}, a Poisson structure is constructed on a Hilbert manifold such that the Poisson bracket of two functions depends on its second-order jets, hence the dual map of the Poisson bracket cannot be factored as a smooth section of $L^2_{\text{skew}}(T^*M,\mathbb R)$. In this article, we assume that the modeling space $E$ of a smoothly paracompact manifold $M$ has the bornological approximation property.
Then, if $\{-,-\}:C^\infty(M)\times C^\infty(M)\to C^\infty(M)$ is a Poisson bracket on $C^\infty(M)$, we show that the dual map ${\widehat{\{\;,\;\}}'}$ factors as a smooth section of the vector bundle $L^2_{\text{skew}}(T^*M,\mathbb R)$. 


\section{The general setting}\nmb0{2}

\subsection{Poisson brackets}\nmb.{2.1}
Let $M$ be a smooth convenient manifold. A Poisson bracket on $M$ is a bounded bilinear mapping 
\begin{align*}
&\{\;,\;\} : C^{\infty}(M) \x C^{\infty}(M) \to C^{\infty}(M) \quad \text{ satisfying  }
\\&
\{f,g\} = -\{g,f\}, \qquad  \{f,gh\} = \{f,g\}h + g\{f,h\}
\\&
\{f,\{g,h\}\} = \{\{f,g\},h\} + \{g,\{f,h\}\}
\end{align*}
for all $f$, $g$, $h \in C^{\infty}(M)$.
The Poisson bracket is said to be of \emph{first order} if $f\mapsto \{f,g\}$ is a linear differential operator of the first order. It is necessary to ask for this since in infinite dimensions there exist derivations of the algebra $C^{\infty}(M)$, which are differential operators of order 2 or 3 (see \cite[28.3]{KrieglMichor97})

We may view this as a bounded linear mapping on the Mackey-completed borno\-logical tensor product (see \cite[5.7]{KrieglMichor97}), which linearizes bounded bilinear mappings and its skew symmetric subspace and its dual mapping
\begin{align*}
&C^{\infty}(M)\bar\otimes_\be C^{\infty}(M) \supset \bigwedge^2 C^{\infty}(M) \xrightarrow{\widehat{\{\;,\;\}}} C^{\infty}(M). 
\\&
\xymatrix{
C^{\infty}(M)' \ar[r]^{\widehat{\{\;,\;\}}'\quad}  & \big(\bigwedge^2 C^{\infty}(M) \big)' \ar@2{-}[r] & L^2_{\text{skew}}(C^{\infty}(M),\mathbb R).
}\end{align*}      

\begin{theorem}\label{main} If $M$ is a smoothly paracompact convenient manifold modeled on a Montel space $E$, such that one of the following holds:
\begin{enumerate}
\item[(a)] $E$ has the bornological approximation property (see \cite[6.6]{KrieglMichor97}), or
\item[(b)] $E$ is a Montel space, and it has the approximation property, i.e.  $E'\otimes E$ is dense in $L_{p}(E,E)$, where $L_{p}(E, E)$ denotes the space of all continuous maps from $E$ to $E$ with the topology of uniform convergence on  pre-compact subsets of $E$.
\end{enumerate}
Then the dual mapping of a Poisson  bracket factors as follows 
\begin{equation*}
\xymatrix{ 
C^{\infty}(M)' \ar[r]^{\widehat{\{\;,\;\}}'\qquad}  & L^2_{\textrm{skew}}(C^{\infty}(M),\mathbb R) \\
M \ar[u]^{\on{ev}} \ar[r]^{P\qquad} & L^2_\textrm{skew}(T^*M,\mathbb R ) \ar[u]_{L^2_{\textrm{skew}}(\on{ev}\o d,\mathbb R)}
}
\end{equation*}
where $\on{ev}_x(f) = f(x)$ and $(L^2_\textrm{}{skew}(\on{ev}\o d,\mathbb R)\bar P_x)(f,g) = \bar P_x(df_x,dg_x)$.
Moreover, the Schouten-Nijenhuis bracket described in \cite{Michor25a} is applicable to $P$, and $[P,P]=0$.
\end{theorem}
\begin{remark}\cite[6.14]{KrieglMichor97}
    The following convenient spaces have the bornological approximation property:
    \begin{itemize}
     \item The space $\Ga(M \leftarrow F)$ of smooth sections of any smooth finite dimensional vector bundle $F \to M$ with separable base $M$.
     \item The space $\Ga(M\leftarrow F)$ of smooth sections with compact support in any finite dimensional vector bundle $F\to M$ with separable base $M$.
    \item  The Fr\'echet space of holomorphic sections  $\Ga_{\mathcal{H}}(M\leftarrow E)$ of a holomorphic vector bundle over a complex Stein manifold $M$. 
    \end{itemize}
\end{remark}
\begin{proof}
We have to show that for any $x\in M$ and $f,g\in C^{\infty}(M)$ the expression $\{f,g\}_x$ depends only on $df_x$ and $dg_x$ in $T^*_xM$. By skew symmetry it suffices to do this for $f$, say. 
\\
Step 1: \emph{If $f=0$ on an open set $U\subset M$, then also $\{f,g\}=0$ on $U$  for any $g$.}
Since $M$ is smoothly paracompact, for any $y\in U$ there exists $h\in C^{\infty}(M)$ such that $h(y)=0$ and $h|_{M\setminus U}=1$. Then $hf=f$ and 
$\{f,g\} = \{hf,g\} = h\{f,g\} + \{h,g\}f$ vanish at $y$.
\\
Step 2: Assume that 
$df_x=0$. Without loss of generality, we can also assume that $f(x)=0$ since for any constant C we have $\{C,g\} = 0$.  
By step 1 we may replace $M$ by a smooth chart $V$ centered at $x$ which is diffeomorphic to an absolutely convex neighborhood of $x=0$ in the modeling convenient vector space $E$.
Let $h\in C^{\infty}(M)$ be a function with support in  $V$  
 with $h=1$ in an open neighborhood $U$ of $0$ in $V$. Then $(h^2f - f)|_U = 0$ thus by step 1 also 
$\{hf-f,g\}|_U = 0$. In $V$ we have 
\begin{align*}
f(y) &= f(0) + df_0(y) + \int_0^1 d^2f(ty)(y,y)(1-t)dt =  \int_0^1 d^2f(ty)(y,y)(1-t)dt
\end{align*} 
Consider the remainder $R_t(y)=d^2f(ty)(y,y)$, for a fixed $t\in [0,1]$ and $y\in U, U{\subset} E $ is open. And \[R(y)=\int_0^1R_t(y) (1-t)dt.\] Then  $R_t, R \in C^\infty(U)$. For each $g\in C^\infty(U), x\in U $, $D_{g,x}:=\{g, -\}(x)$ is a point derivation, which is a linear function on $C^\infty(U)$. Since linear functions can be passed under the integral, we have, \[ D_x(R)=\int_0^1D_x(R_t) (1-t)dt .\]
Now its enough to show that for each $t\in[0,1]$, $D_x(R_t)=0$. Without loss of generality, we assume $x=0$.  

In case \thetag{a}, since $E$ has the bornological approximation property, $id_E\in L(E,E)$ can be approximated by elements in $E'\otimes E$, there exists a net of finite-dimensional linear operators $L_i \in E'\otimes E$ converges to $id_E  \in L(E,E)$ for the topology of uniform convergence on  bounded subsets. In case \thetag{b}, when $E$ has the usual approximation property, we can find a net $L_i \in E'\otimes E$ which converges to $id_E  \in L(E,E)$ only uniformy on compact subsets. But since $E$ is also a Montel space, this convergence is the same as the uniform convergence on bounded subsets of $E$. 

Let \[R_{t,i}(y):=d^2f(ty)(L_iy,y).\]
For fixed $t\in [0,1]$ we have $R_{t,i}\in C^\infty(U)$;  we want to show that $R_{t,i} \to R_t$ in the convenient topology of  $C^\infty(U)$.

 Let $c\in C^\infty(\mathbb R, U)$ be a smooth curve, then 
$R_{t}\circ c$ and $R_{t,i}\circ c\in C^\infty(\mathbb R)$; we need to show that $R_{t,i}\circ c$ and its all derivatives converge uniformly on compact intervals in $R$ to $R_{t,i}\circ c$.

Let $[a,b]\subset \mathbb R$ be a compact interval. The map $\mathbb R\to E'$ defined by: $y\to d^2(f\circ c)(ty)(-,y)$ is a smooth curve and is thus continuous (by \cite[1.8]{KrieglMichor97}, e.g.),
so the image $C$ of $[a,b]$ under this map is compact, hence the image is weakly bounded, but since a convenient space $E$ is quasibarrelled, 
the image is equicontinuous in $E'$; see \cite[III 4.2]{schaefertopological}. This means that \[ |d^2(f\circ c)(ty)(v,y)| < \epsilon \]
 for $y\in C, v\in V$, where $V$ is an open $0$ neighborhood of $E$. Hence, \[
 |R_{t,i}\circ c (y)-R_{t}\circ c (y)| <\epsilon \]
 for $y\in [a,b].$ 
 Therefore, $R_{t,i}\circ c\to R_t\circ c$ uniformly on compact intervals in $\mathbb R$. We can use the same argument for derivatives of $R_{t,i}\circ c $ and $R_{t}$, and conclude that $R_{t,i}\to R_{t}$ in $C^\infty(U)$. Now, since $L_i\in E'\otimes E$ is a  finite rank operator, $L_i=\sum_{j=1}^{j=i}e'_{ij}\otimes e_{ij}$, for $e_{ij}' \in E'$ and $e_{ij}\in E$. Then \[D_0(R_{t,i})(y)=\sum_{j=1}^{i} D_0(d^2f(y)(e'_{ji}(y)e_{ij}, y)=\sum_{j=1}^{i} D_0((e'_{ji}(y))(d^2f(y)e_{ij}, y)).\]
Now, let us define $h_{ij}(y):=e'_{ij}(y), g_{ij}(y):=d^2f(y)(e_{ij}, y) $. Then, $h_{ij}, g_{ij}$ are smooth functions on a suitable $0$-neighborhood $V\subset E.$ Now, 
\[D_0(R_{t,i})=\sum_{j=1}^{i}D_{0}(h_{ij})g_{ij}(0)+D_{0}(g_{ij})h_{ij}(0)=0.\] As $R_{t,i} \to R_t$ in $C^\infty(U)$, and $D_0$ is a bounded linear functional on $C^\infty(U)$, we have $D_{0}(R_{t})=\underset {i\to \infty}{\lim} D_{0}(R_{t,i})=0.$ 
\end{proof}

\begin{corollary}
    A Montel space with the approximation property has the bornological approximation property.
\end{corollary}
\begin{remark}
    Examples of Montel space with the approximation property are spaces of smooth sections of finite dimensional vector bundles with compact base, or even noncompact base when one considers inductive limits of all subspaces of sections with compact support; spaces of holomorphic sections of complex vector bundles; but in general not their dual spaces.
\end{remark}
\begin{corollary}
    If $M$ is a smoothly paracompact manifold modeled on a nuclear Fr\'echet (or a Silva space) then the dual mapping of a Poisson  bracket factors as a smooth section of $L^2_{skew}(T^*M)\cong \wedge^2TM$ . Where $\wedge^2TG$ denotes the smooth vector bundle with the fiber $\wedge^2T_xM$ over $x\in M$; the skew-symmetric projective tensor product of $T_xM$ with itself.

\end{corollary}
\begin{remark}
   Note that if  $M$ is modeled on a nuclear Fréchet space (see \cite{thomas95}) or on a nuclear (LF) space (see \cite{michor1983manifolds}), similar arguments are used to show the topological isomorphism between the space of continuous derivations on $C^\infty(M)$ and the space of smooth vector fields $\mathfrak{X}(M)$.
\end{remark}

\

\section{Examples}\nmb0{3}

\subsection{The Poisson structure on the dual of a Lie algebra}
\nmb.{3.1}
Let $\mathfrak g$ be a convenient Lie algebra; i.e., $\mathfrak g$ is a convenient vector space 
together with a bounded Lie bracket $[\;,\;]:\mathfrak g\x \mathfrak g\to \mathfrak g$.
By the universal property of the bornological tensor product, the Lie bracket extends to a linear 
bounded mapping 
$$\xymatrix{
\mathfrak g\x \mathfrak g  \ar[r]^{[\;,\;]} \ar[d] & \mathfrak g
\\
\bigwedge\nolimits^2\mathfrak g \ar[ur]_{\widetilde{[\;,\;]}}
}$$
The dual mapping 
$$
P:= \widetilde{[\;,\;]}^*:\mathfrak g^*\to  (\bigwedge\nolimits^2\mathfrak g)^* = 
L^2_{\text{skew}}(\mathfrak g;\mathbb R)
$$
is the candidate for the Lie-Poisson structure on the dual of the Lie algebra.
It induces a Lie bracket on the space 
 $$C^{\infty}_P(\mathfrak g^*)=\{f\in C^{\infty}(\mathfrak g^*): df \text{ factors as } i_{\mathfrak g}\o\widetilde{df} \text{ for  } \widetilde{df}\in C^{\infty}(\mathfrak g^*,\mathfrak g)\}$$
with the convenient vector space structure induced by the linear mappings 
$$ C^{\infty}_P(\mathfrak g^*) \xrightarrow{\text{incl}} C^{\infty}(\mathfrak g^*,\mathbb R),\quad  
 C^{\infty}_P(\mathfrak g^*) \xrightarrow{d} C^{\infty}(\mathfrak g^*,\mathfrak g)\,.$$
 Similarly to \cite[48.6]{KrieglMichor97} we have the following result:
 
\begin{lemma} For $f \in C^\infty(\mathfrak g^*,\mathbb R)$
we have \thetag{1} $\Leftrightarrow$ \thetag{2} $\Rightarrow$ \thetag{3}, where:
\begin{enumerate}
\item $df:\mathfrak g^* \to  \mathfrak g^{**}$ actually factors as a smooth mapping 
$\mathfrak g^* \xrightarrow{df}  \mathfrak g \xrightarrow{i_{\mathfrak g}} \mathfrak g^{**}$. 
\item Each iterated derivative $d^kf_x \in L^k_{\text{sym}}(\mathfrak g^*,\mathbb R)$ has the property that 
	$$d^kf_x(\quad,y_2,\dots,y_k) \in \mathfrak g \text{ for all  }x,y_2,\dots,y_k \in \mathfrak g^*,$$
	and that the mapping 
\begin{align*}
\prod^k\mathfrak g^* \to  \mathfrak g,\quad\text{ given by } \quad
(x,y_{2},\dots,y_{k})\mapsto df_x(\quad,y_{2},\dots,y_{k})\,,
\end{align*}
	is smooth. By the symmetry of higher derivatives, this is then true for all entries of $d^kf_x$, for all $x$.  
\item $f$ has a smooth $P$-gradient 
        $\on{grad}^P f = P\o df\in \X(\mathfrak g^*) = C^\infty(\mathfrak g^*,\mathfrak g^*)$ which 
        is given by $P_x(df_x,y) = \langle y,\on{grad}^Pf(x)\rangle$. 
\end{enumerate}
\end{lemma}

This looks like the usual finite dimensional version if $\mathfrak g$ is reflexive. 

For non-reflexive $\mathfrak g$ let $i=i_{\mathfrak g}: \mathfrak g \to \mathfrak g^{**}$ be canonical embedding into the bidual.
Then we may ask for a linear lift $\tilde P$ fitting into the following diagram
$$
\xymatrix{
\mathfrak g^* \ar[r]^{P:= \widetilde{[\;,\;]}^*} \ar@{-->}[dr]_{\tilde P} 
& (\bigwedge\nolimits^2\mathfrak g)^* \ar@{=}[r] & L^2_{\text{skew}}(\mathfrak g;\mathbb R)
\\
& (\bigwedge\nolimits^2\mathfrak g^{**})^* \ar[u]_{\bigwedge\nolimits^2(i_{\mathfrak g})^*} \ar@{=}[r]
& L^2_{\text{skew}}(\mathfrak g^{**};\mathbb R) \ar[u]_{L^2_{\text{skew}}(i_{\mathfrak g};\mathbb R)}
}
$$
If $\mathfrak g$ has predual $\mathfrak g^p$ with $(\mathfrak g^p)* = \mathfrak g$ as a convenient vector space, we get such a natural lift: 
Let $i_{\mathfrak g^p}:\mathfrak g^p \to (\mathfrak g^p)^{**}=\mathfrak g^*$ be the canonical injection,
$(i_{\mathfrak g^p})^*: \mathfrak g^{**} \to \mathfrak g$, so that
$$
\xymatrix{
\mathfrak g^* \ar[r]^{P:= \widetilde{[\;,\;]}^*\quad} \ar@{-->}[dr]_{\tilde P}  
& (\bigwedge\nolimits^2\mathfrak g)^* \ar@{=}[r]  \ar[d]^{(\bigwedge\nolimits^2(i_{\mathfrak g^p})^*)^*} 
& L^2_{\text{skew}}(\mathfrak g;\mathbb R) \ar[d]^{L^2_{\text{skew}}((i_{\mathfrak g^p})^*,\mathbb R)}
\\
& (\bigwedge\nolimits^2\mathfrak g^{**})^* 
\ar@{=}[r]
& L^2_{\text{skew}}(\mathfrak g^{**};\mathbb R) 
}
$$

\subsection{Example: Lie algebras of vector fields on $\mathbb R^n$}\nmb.{3.2}
We first consider the Lie algebra $\X_c(\mathbb R^n)$ of vector fields with compact support, a 
nuclear (LF)-space. By \cite[Theorem 51.6]{Treves67} we have 
$L(\X_c(\mathbb R^n),\X_c(\mathbb R^n)')= 
\mathcal D'(\mathbb  R^n)^n\widehat{\otimes}_{\pi}\mathcal D'(\mathbb R^n)^n$.
Thus, we get the Lie Poisson structure:
\begin{align*}
P=-\widetilde{[\;,\;]}^*: \mathcal D'(\mathbb R^n)^n =\Ga_{\mathcal D'}(T^*\mathbb R^n)\to 
(\bigwedge\nolimits^2 \X_c(\mathbb R^n))^* = \bigwedge\nolimits^2 \mathcal D'(\mathbb R^n)^n.
\end{align*}
It fits the setup of section \nmb!{1}, thanks to nuclearity, $\X_c(\mathbb R^n)$ has the approximation property (see \cite[22.2, Corollary 2]{jarchow2012locally} ), and since $\X_c(\mathbb R^n)$ is reflective $D'(\mathbb R^n)=\X_c(\mathbb R^n)'$ also has the approximation property by \cite[Lemma 6.1]{KrieglMichor97}.
We compute $P(A)$ for a distributional 1-form (a 1-co-current) 
$A=\sum_j A_j.dx^j\in \Ga_{\mathcal D'}(T^*\mathbb R^n)$. It suffices to compute its action on 
$X\otimes Y$ for $X=\sum_i X^i\p_i,Y=\sum_iY^i\p_i\in\X_c(\mathbb R^n)$. We shall make use of the 
description of distributions on $\mathbb R^n\x \mathbb R^n$  which are supported on the diagonal 
given in \cite[Theorem 5.2.3]{Hormander83} which is based on \cite[Theorem 2.3.5]{Hormander83}. We 
will give enough details to make the use of these results clear, and 
we will use the diagonal mapping $\on{diag}:\mathbb R^n\to \mathbb R^n\x\mathbb R^n$.
\begin{align*}
\langle P(A)&,X\otimes Y \rangle = \langle A, -[X,Y] \rangle 
= \sum_{i,j} \langle A_j, Y^i\p_i X^j - X^i\p_i Y^j \rangle
\\&
= \sum_{i,j} \big(\langle A_i, (\p_j X^i)Y^j\rangle - \langle A_j, X^i(\p_i Y^j) \rangle
\\&
= \sum_{i,j} \big(\langle A_i, (\p_j X^i)Y^j\rangle - \langle A_j, X^i(\p_i Y^j) \rangle
\\&
= \sum_{i,j} \Big(\langle A_i, (\tfrac{\p}{\p x_j} X^i(x)Y^j(y))\o\on{diag}\rangle 
                  - \langle A_j, (\tfrac{\p}{\p y^i}X^i(x)Y^j(y))\o\on{diag} \rangle\Big)
\\&
= \sum_{i,j} \Big(\langle A_i, \on{diag}^*(\tfrac{\p}{\p x_j} X^i(x)Y^j(y))\rangle 
                  - \langle A_j, \on{diag}^*(\tfrac{\p}{\p y^i}X^i(x)Y^j(y)) \rangle\Big)
\\&
= \sum_{i,j} \Big\langle \big(\on{diag}_* A_i\big) \tfrac{\p}{\p x_j} 
	-\big(\on{diag}_* A_j)\tfrac{\p}{\p y^i}, X^i(x)Y^j(y)) \Big\rangle
\end{align*}
so that
\begin{equation}\boxed{\quad
P\big(\sum\nolimits_k A_k dx^k\big) = \sum\nolimits_{i,j} \big((\on{diag}_* A_i) \tfrac{\p}{\p x_j} 
	-(\on{diag}_* A_j)\tfrac{\p}{\p y^i}\big)\cdot dx^i\otimes dy^j
\quad}\label{eq:1}\end{equation}
Analogous results also hold for the Lie algebras $\X_{\mathcal S}(\mathbb R^n)$ of rapidly 
decreasing smooth vector fields, and for the Lie algebra $\X(\mathbb R^n)$ of all smooth vector 
fields, since they are also nuclear. For other Lie algebras like $\X_{H^\infty}(\mathbb R^n)$ or 
$\X_{\mathcal B}(\mathbb R^n)$ (in the notation from \cite{MichorMumford2013z}) formula \eqref{eq:1} for $P$ 
still holds, but the image is not in the space of summable skew multi vector fields, since the 
spaces are not nuclear. The same holds for the Lie algebras based on Denjoy-Carleman 
ultradifferentiable functions on $\mathbb R^n$ of Roumieu and Beurling type discussed in \cite{KrieglMichorRainer14}, 
depending on whether they are nuclear or not. For the non-nuclear ones, the Poisson structure is, 
however, a bounded skew multi vector field on the dual of the Lie algebra. 

All the above also holds if we replace $\mathbb R^n$ by an open subset therein, and formula \eqref{eq:1}; this paves the way 
to carry this over to manifolds. We summarize this as:

\begin{proposition}
Let $N$ be a finite-dimensional smooth manifold. Then the Poisson structure on the dual $\X(N)'$ of 
the Lie algebra $\X(N)$ of all smooth vector fields on $N$,
given as 
$$
P=-\widetilde{[\;,\;]}^*: \X(N)' =\Ga_{\mathcal D'}(T^*N)\to 
(\bigwedge\nolimits^2 \X_c(N))' = \bigwedge\nolimits^2 \X_c(N)'
$$
is given by formula \eqref{eq:1} on each open chart of $N$.

The same holds for all other suitable Lie algebras of vector fields: For example, for the Lie algebra $\X_\mu(M)$ of divergence free vector fields on a (say) compact oriented manifold $(M,\mu)$ with positive smooth volume form $\mu$, the dual space is 
$$\X_\mu(M)' = \Ga_{\mathcal D'}(T^*M)/ d(\mathcal D'(M)),$$
and formula \eqref{eq:1} factors to the quotient. 
\end{proposition} 

\section*{Acknowledgments}
The second author was supported by the German Research Foundation (Deutsche Foschungsgemeinschaft), SFB-TRR 358/1 2023 – 491392403.


\end{document}